\documentclass[12pt]{article}
\usepackage{amsfonts,amssymb}

\usepackage{amsmath,bbm}
\usepackage{theorem}          
\usepackage{latexsym}
\usepackage{float}                   
\usepackage{times}
\usepackage{cite}

\newtheorem{tht}{Theorem}[section]

\newtheorem{thl}[tht]{Lemma}

\newcommand{\anf}{\raisebox{0.2ex}{\scriptsize$\triangleleft$}}
\newcommand{\ang}{\raisebox{0.2ex}{\scriptsize$\triangleright$}}

\newcommand{\sn}{\smallskip\noindent}

\newcommand{\lti}{\,{\scriptstyle\ltimes}\,} 
\newcommand{\rti}{\,{\scriptstyle\rtimes}\,} %

\newcommand{\cg}{\mathcal{G}} 
 
\newcommand{\cO}{{\mathcal{O}}}
\newcommand{\cK}{{\mathcal{K}}}
\newcommand{\cU}{{\mathcal{U}}}   
\newcommand{\cT}{{\mathcal{T}}}

\newcommand{\cX}{{\mathcal{X}}}

\newcommand{\cS}{\mathcal{S}} 
\newcommand{\cL}{\mathcal{L}}

\newcommand{\cJ}{\mathcal{J}}

\newcommand{\dC}{\mathbb{C}}
\newcommand{\dN}{\mathbb{N}}

\newcommand{\im}{\mathrm{i}}
\newcommand{\Lin}{{\mathrm{Lin}}}

\newcommand{\dd}{\mathrm{d}}

\newcommand{\TT}{\mathrm{t}}

\newcommand{\tr}{\mathrm{Tr}}
\newcommand{\Sq}{\mathrm{S}_q^2}

\newcommand{\SU}{\mathrm{SU}_q(2)}
\newcommand{\su}{\mathrm{su}_2}

\newcommand{\RT}{\mathrm{Re}}

\newcommand{\comp }{\mathbb{C}}

\newcommand{\dif }{\mathrm{d}}

\newcommand{\DCS }{\Gamma }

\newcommand{\lid }{\mathcal{L}}
\newcommand{\OSUq }[1]{\mathcal{O}(\mathrm{SU}_q(2))}
\newcommand{\ot }{\otimes }
\newcommand{\otp }{\ot _{\podl }}

\newcommand{\pair }[2]{\langle #1,#2 \rangle }
\newcommand{\pinv }{\omega _{\mathrm{R}}}
\newcommand{\podl }{\mathcal{O}(\mathrm{S}^2_q)}

\newcommand{\tpair }[2]{\langle \!\langle #1,#2 \rangle \!\rangle }
\newcommand{\vep }{\varepsilon }

\begin{document}

\date{\small{Fakult\"at f\"ur Mathematik und 
Informatik\\ Universit\"at Leipzig, 
Augustusplatz 10, 04109 Leipzig, Germany\\ 
E-mail: schmuedg@mathematik.uni-leipzig.de / 
wagner@mathematik.uni-leipzig.de}
}

\title{Dirac operator and a twisted cyclic cocycle on the standard 
Podle\'s quantum sphere}

\author{Konrad Schm\"udgen and Elmar Wagner}

\maketitle

\renewcommand{\theenumi}{\roman{enumi}}
\begin{abstract}
A Dirac operator $D$ on the standard Podle\'s sphere $\Sq $ is defined and 
investigated. It yields a real spectral triple such that $|D|^{-z}$ is of 
trace class for $\RT\, z>0$. Commutators with 
the Dirac operator give the distinguished 2-dimensional covariant differential 
calculus on $\Sq $. The twisted cyclic cocycle associated with the 
volume form of the differential calculus is expressed by means of the 
Dirac operator.
\end{abstract}
{\it Keywords:} Quantum spheres, differential calculus, spectral triples, 
twisted cyclic   
cocycles\\
{\it Mathematics Subject Classifications (2000):} 
17B37, 46L87, 58B37, 81R50
%

%
\setcounter{section}{-1}
\section{Introduction}
Spectral triples, cyclic cocycles and Dirac operators are basic concepts 
of Alain Connes' non-commutative geometry \cite{C}. In general, covariant 
differential calculi on quantum groups cannot be described by 
spectral triples.
For instance, the commonly used covariant differential calculi on 
the quantum group $\SU$ cannot even be given by bounded commutators 
\cite{S1}. For the quantum group $\SU$, an equivariant spectral triple 
has been constructed in \cite{CP} and the corresponding non-commutative 
geometry has been studied by A. Connes in \cite{C3}. 
In this paper we show (among other things) 
that the distinguished 
2-dimensional covariant differential calculus of the standard 
Podle\'s quantum sphere $\Sq $ can be given by a spectral triple.

The starting point of our construction is the representation of the left 
crossed product algebra $\cO(\Sq )\rti \cU_q(\su )$ of the 
Hopf algebra $\cU_q(\su )$ and the coordinate algebra $\cO(\Sq )$ of the 
standard Podle\'s sphere on the subspace $V^+ \oplus V^-$ of $\cO(\SU)$, 
where $V^+=\cO(\Sq ) c+ \cO(\Sq )d$ and $V^-=\cO(\Sq )a+\cO(\Sq )b$.
Then the right actions $R_E$ and $R_F$ of the generators 
$E,F\in\cU_q(\su )$ on 
$\cO(\SU)$ map $V^+$ into $V^-$ and $V^-$ into $V^+$, respectively. 
The Dirac operator
\begin{equation}
D=\left( \begin{matrix} 0 &R_F\\ R_E & 0 \end{matrix}\right)
\end{equation}
commutes with the left action of $\cU_q(\su )$ on $V^+\oplus V^-$. 
Since $D$ has discrete spectrum of eigenvalues $[l+1]_q$, $l\in\dN_0$, 
with multiplicities $2l+1$, and since commutators $[D,x]$, $x\in \cO(\Sq )$, 
are bounded, the Dirac operator $D$ leads to a spectral triple such that 
$|D|^{-z}$ is of trace class for all $z\in \dC$, $\RT\, z>0$. 

Note that our Dirac operator $D$ is unitarily equivalent to the Dirac 
operator constructed in \cite{DS}, where also a real structure 
was obtained. However, we show that 
the Dirac operator fits into the framework of covariant differential 
calculus. That is, we prove that the Dirac operator gives a commutator 
representation  $\dd x\cong \im [D,x]$ of the distinguished 
2-dimensional covariant 
first order differential calculus on $\cO(\Sq )$. The associated higher 
order calculus has up to multiples a unique invariant volume 2-form $w$. 
Let $h$ denote the $\cU_q(\su )$-in\-var\-i\-ant state on $\cO(\Sq )$. We show 
that the associated twisted cyclic 2-cocycle $\tau_{w,h}$ represents a 
non-trivial cohomology class in the twisted cyclic cohomology 
and prove that
\begin{align*}
\tau_{\omega,h} (x_0, x_1, x_2)
&=h(x_0(R_F (x_1) R_E (x_2) - q^2 R_E (x_1) R_F(x_2)))\\
&=(q-q^{-1})^{-1} (\log q)~ \underset{z=2}{{\rm res}}\, 
\tr_\cK\, \gamma_q K^2 |D|^{-z} \,x_0[D,x_1][D,x_2] 
\end{align*}
for $x_0, x_1,x_2\in\cO(\Sq )$, where
\begin{equation*}
\gamma_q= \left(\begin{matrix} 1 &0\\ 0 &-q^2 \end{matrix}\right).
\end{equation*}
This paper is organized as follows. In Section \ref{S1} 
twisted cyclic cocycles 
associated with volume forms of differential calculi are defined. 
A number of preliminary facts on the Podle\'s sphere $\Sq $ and the 
quantum group $\SU$ are collected in Section \ref{S2}. 
The Dirac operator and its properties are developed in Section \ref{S3}, 
while the twisted cyclic cocycle is studied in Section \ref{S4}. 
The form of the cocycle $\tau_{\omega,h}$ was kindly 
communicated to us by I. Heckenberger. His proof of the corresponding 
result will be given in the appendix.

Throughout this paper, $q$ is a real number such that $0<q<1$ and 
$\im$ denotes the complex unit. 
We abbreviate 
$\lambda := q - q^{-1}$ and $[n]_q := (q^n - q^{-n})/(q-q^{-1})$.

\section{Twisted cyclic cocycles}                                \label{S1}
Let $\cX$ be a unital algebra and $\sigma$ an algebra automorphism of $\cX$. 
We shall use the following notations on the twisted cyclic cohomology 
resp.\ twisted cyclic homology (see e.g.\ \cite{KMT} and \cite{KR}) of the 
pair $(\cX,\sigma)$.

Let $\varphi$ be an $(n+1)$-linear form on $\cX$ and let  
$\eta = \sum_k x^k_0\otimes {\cdots}\otimes x^k_n \in \cX^{\otimes n+1}$, 
where $x_j^k\in\cX$. The $\sigma$-twisted coboundary resp.\ boundary 
operator $b_\sigma$ and the $\sigma$-twisted cyclicity operator 
$\lambda_\sigma$ are defined by
\begin{align}
(b_\sigma\varphi)(x_0,{\dots},x_{n+1})                    \nonumber
&=\sum^{n}_{j=0} (-1)^j \varphi(x_0,{\dots},x_j x_{j+1},{\dots}, x_{n+1})\\
& \qquad \qquad\quad\quad                            \label{b1}
+(-1)^{n+1} \varphi (\sigma (x_{n+1}) x_0,x_1,{\dots},x_{n}),\\
(\lambda_\sigma\varphi)(x_0,{\dots},x_{n})           \label{lambda}
&=(-1)^{n} \varphi(\sigma(x_n),x_0,{\dots}, x_{n-1}),\\
b_\sigma\eta                                         \nonumber
&=\sum_k \sum^{n-1}_{j=0} (-1)^j x^k_0\otimes{\cdots}
\otimes x^k_j x^k_{j+1}\otimes {\cdots}\otimes x^k_n\\
& \qquad \qquad\quad\quad 
+(-1)^n \sigma (x^k_n) x^k_0 \otimes x^k_1 
\otimes{\cdots}\otimes x^k_{n-1},\\
\lambda_\sigma\eta 
&=\sum_k (-1)^{n}\sigma(x^k_n)\otimes x^k_0 \otimes {\cdots}\otimes x^k_{n-1}.
\end{align}
An $(n+1)$-form $\varphi$ is called a 
{\it $\sigma$-twisted cyclic n-cocycle} if 
$b_\sigma \varphi=0$ and $\lambda_\sigma\varphi=\varphi$. 
A $\sigma$-twisted cyclic $n$-cocycle $\varphi$  
is called  {\it non-trivial} if 
there is no $n$-form $\psi$ such that $\varphi=b_\sigma \psi$ and 
$\lambda_\sigma \psi=\psi$. Similarly, $\eta$ is a 
$\sigma$-twisted cyclic $n$-cycle if $b_\sigma\eta=0$ and 
$\lambda_\sigma\eta=\eta$. Cycles 
$\eta = \sum_k x^k_0\otimes {\cdots}\otimes x^k_n \in \cX^{\otimes n+1}$
can be paired with $(n+1)$-forms $\varphi$ by setting 
$\varphi(\eta):=\sum_k \varphi(x^k_0,{\ldots}, x^k_n)$.

Let $\Gamma={\bigoplus^\infty_{k=0}} \Gamma^{\wedge k}$ be 
a differential calculus on $\cX$ with differentiation 
$\dd:\Gamma^{\wedge k}\rightarrow \Gamma^{\wedge(k+1)}$
(see e.g. \cite[Part IV]{KS}). 
Suppose that $\Gamma^{\wedge n}$ is a free $\cX$-module generated by 
$\omega$ such that  $\omega x = x\omega$ 
for all $x\in\cX$. 
For $\eta\in\Gamma^{\wedge n}$, let $\pi(\eta)\in\cX$ denote the unique 
element such that $\eta=\pi(\eta)\omega$. 
Assume that $\sigma$ is an algebra 
automorphism of $\cX$ and $h$ is a linear functional on $\cX$ such that,
for all $x,y\in\cX$ and $\omega_{n-1}\in\Gamma^{\wedge(n-1)}$,
$$
h(xy) = h(\sigma(y)x), \quad 
h(\pi(\dd \omega_{n-1}))=0.
$$
Then it is easy to check 
by using repeatedly the Leibniz rule that
$$
\tau_{\omega,h} (x_0,{\dots},x_n) 
:= h(\pi (x_0 \dd x_1 \wedge{\dots}\wedge \dd x_n)),\quad  
x_0,{\dots}, x_n \in \cX,
$$
defines a $\sigma$-twisted cyclic $n$-cocycle $\tau_{\omega,h}$ on $\cX$.
%
%
\section{Preliminaries on the quantum group ${\bf \SU}$ and the 
Podle\'s sphere ${\bf \Sq }$}                                   \label{S2}
%
First we collect some definitions, 
facts and notations  used in what follows
(see e.g.\ \cite{KS}).  
For brevity, we employ the Sweedler notation 
$\Delta(x)={x_{(1)}\otimes x_{(2)}}$ for the  
comultiplication $\Delta(x)$.

Let $\cO(\SU)$ denote the coordinate Hopf $\ast$-al\-ge\-bra of 
the quantum $\mathrm{SU}(2)$ group (see e.g.\ \cite[Chapter 4]{KS}). As usual, 
the generators of $\cO(\SU)$ are denoted by $a$, $b$, $c$, $d$. 
The Hopf $\ast$-al\-ge\-bra $\cU_q(\su )$ has four generators 
$E$, $F$, $K$, $K^{-1}$ with defining relations
$$
KK^{-1}=K^{-1} K=1,\ KE=q EK,\ 
FK= qKF,\ EF-FE=\lambda^{-1} (K^2-K^{-2}),
$$
involution $E^\ast=F$, $K^\ast=K$, comultiplication 
$$\Delta (E)=E\otimes K+K^{-1} \otimes E, \ \ 
\Delta(F) = F\otimes K+K^{-1} \otimes F,\ \  \Delta (K)=K\otimes K,
$$
counit $\varepsilon (E) = \varepsilon (F)= \varepsilon (K-1)=0$ and 
antipode  $S(K)=K^{-1}$, $S(E)=-qE$, $S(F)=-q^{-1}F$. 
There is a dual pairing $\langle\cdot,\cdot\rangle$ of the 
Hopf $\ast$-al\-ge\-bras $\cU_q(\su )$ and $\cO(\SU)$ given 
on the generators by the values
\begin{equation}\label{pair}
\langle K^{\pm 1}, d\rangle = \langle K^{\mp 1},a\rangle = q^{\pm 1/2},\quad 
\langle E,c\rangle=\langle F,b\rangle =1
\end{equation}
and zero otherwise. Therefore, $\cO(\SU)$ is a left and right 
$\cU_q(\su )$-module $\ast$-al\-ge\-bra with left action $\ang$ and right 
action $\anf$ defined by
\begin{equation}\label{act}
f\ang x=\langle f, x_{(2)}\rangle x_{(1)},\ \  
x\anf f=\langle f, x_{(1)}\rangle x_{(2)},\quad  x\in \cO(\SU), \ 
f\in \cU_q(\su ).
\end{equation} 
The actions satisfy 
\begin{align}                                       \label{actstar}
&(f\ang x)^\ast=S(f)^\ast\ang x^\ast,\quad 
(x\anf f)^\ast=x^\ast\anf S(f)^\ast,\\
&f\ang xy=(f_{(1)}\ang x) (f_{(2)}\ang y),\quad      \label{actmod}
xy\anf f=(x\anf f_{(1)})(y\anf f_{(2)}).
\end{align} 
On generators, we have 
\begin{align}\label{leftact1}
&E\ang a=b,\, E\ang c=d,\, E\ang b=E\ang d=0,\, F\ang b=a,\, F\ang d=c,\, 
F\ang a=F\ang c=0,\\
                                                      \label{leftact2}
&K\ang a = q^{-1/2} a,\  K\ang b=q^{1/2} b,\  K\ang c=q^{-1/2}c,\  
K\ang d = q^{1/2} d;\\[6pt]
\label{rightact1}
&c\anf E=a,\, d\anf E=b,\, a\anf E=b\anf E=0,\, a\anf F=c,\, b\anf F=d,\,  
c\anf F=d\anf F=0,\\
                                                     \label{rightact2}
&a\anf K = q^{-1/2} a,\  b\anf K=q^{-1/2} b,\  c\anf K=q^{1/2}c,\  
d\anf K = q^{1/2} d.
\end{align}

Let $h$ denote the Haar state of $\cO(\SU)$ and 
let $\cL^2 (\SU)$ be the Hilbert space completion of $\cO(\SU)$ 
equipped with the scalar product $(x,y)=h(y^\ast x)$, $x, y \in \cO(\SU)$. 
For $x \in \cO(\SU)$ and  $f\in \cU_q(\su )$, set
\begin{equation}                                          
R_f(x)= x\anf S^{-1}(f).                          \label{R}
\end{equation}       
Using (\ref{actstar}) and the $\cU_q(\su )$-invariance of $h$, 
we compute
\begin{align*}
\big( x,R_f(y)\big)&=h((y\anf S^{-1}(f))^\ast x)=h((y^\ast\anf f^\ast)x)
=\varepsilon (f_{(1)}^\ast)h((y^\ast\anf f_{(2)}^\ast)x)\\
&=h((y^\ast\anf f_{(3)}^\ast S^{-1}(f_{(2)}^\ast))
(x\anf S^{-1}(f_{(1)}^\ast)))=\big( R_{f^\ast}(x),y\big).
\end{align*}
Hence the mapping $R$ given by $f\mapsto R_f$ 
is a $\ast$-re\-pre\-sen\-tation of $\cU_q(\su )$. 
Likewise,  the left action $\ang$ defines a 
$\ast$-re\-pre\-sen\-tation of $\cU_q(\su )$ on $\cO(\SU)$.

Let $t^l_{jk}$, $l\in\frac{1}{2}\dN_0$,  $j,k=-l,-l+1,{\dots},l$, denote 
the matrix elements of finite 
dimensional unitary corepresentations of $\cO(\SU)$. 
By the Peter-Weyl theorem, the elements 
$v^l_{jk}:=[2l+1]^{1/2}_q q^jt^l_{jk}$ form an 
orthonormal vector space basis of $\cO(\SU)$. 
For our considerations 
below, we shall need the elements $t^{l+1/2}_{\pm 1/2,j}$, $l\in\dN_0$, 
$j=-(l+1/2),{\dots},l+1/2$. 
They can be expressed explicitely by the following formulas 
(see e.g.\ \cite[p.\ 109]{KS}):
\begin{align}                                            \label{t1}
t^{l+1/2}_{1/2,j} &= M^l_j p_{l-j+1/2} (\zeta; q^{2(1/2-j)}, 
q^{-2(1/2+j)}) b^{j-1/2} d^{j+1/2}, &  j&>0,\\
t^{l+1/2}_{1/2,j} 
&= N^l_{-j} a^{-j-1/2} c^{-j+1/2} p_{l+j+1/2}
(\zeta; q^{-2(1/2-j)}, q^{2(1/2+j)}), & j&<0,\\
t^{l+1/2}_{-1/2,j} &= N^l_j p_{l-j+1/2} 
(\zeta; q^{-2(1/2+j)}, q^{2(1/2-j)}) b^{j+1/2} d^{j-1/2},& j&>0,\\
t^{l+1/2}_{-1/2,j} 
&= M^l_{-j} a^{-j+1/2} c^{-j-1/2}  p_{l+j+1/2} 
(\zeta; q^{2(1/2+j)}, q^{-2(1/2-j)}), &   j&<0,          \label{t4}
\end{align}
where $p_k(\cdot;\cdot,\cdot)$ are the little $q$-Jacobi polynomials, 
$\zeta=-qbc$ and $M^l_j,N^l_j$ are positive constants.
On the vector space 
$W^l_k := {\rm Lin} \{v^l_{jk}\,;\, j=-l,{\dots},l\}$,  
the $\ast$-re\-pre\-sen\-tation $R$ becomes
a spin $l$ representation. 
Let $\alpha^l_j:=([l-j]_q[l+j+1]_q)^{1/2}$. Then  
\begin{equation}                                          \label{tefkact}
R_E(v^l_{jk}) =-\alpha_{j-1}^l v^l_{j-1,k},\quad
R_F(v^l_{jk})=-\alpha_{j}^l v^l_{j+1,k}, \quad
R_K(v^l_{jk})=q^{-j} v^l_{jk}.
\end{equation}
For the left action $\ang$ of $\cU_q(\su )$ on $\cO(\SU)$, one obtains 
\begin{equation}                                            \label{EFK}
E\ang v^{l}_{jk} = \alpha_k^{l} v^{l}_{j, k+1},\quad 
F \ang v^{l}_{j k} = \alpha^{l}_{k-1} v^{l}_{j,k-1},\quad 
K \ang v^{l}_{jk}= q^k  v^{l}_{jk}.
\end{equation}

The coordinate $\ast$-al\-ge\-bra $\cO(\Sq )$  of the standard 
Podle\'s sphere is the unital $\ast$-sub\-al\-ge\-bra of $\cO(\SU)$  
generated by the elements $A:=-q^{-1} bc$, $B:= ac$, 
$B^\ast :=-db$. One can also define $\cO(\Sq )$  as the abstract 
unital $\ast$-al\-ge\-bra with three gen\-e\-ra\-tors 
$A=A^\ast$, $B$, $B^\ast$ and 
defining relations
\begin{equation}\label{algrel}
BA=q^2AB,\ \   AB^\ast = q^2 B^\ast A,\ \   B^\ast B=A-A^2, \ \  
BB^\ast = q^2 A-q^4 A^2.
\end{equation}
It is a well-known fact and easy to verify that $\cO(\Sq )$ is the set 
of all elements of $\cO(\SU)$  which are invariant under the right 
action of the generator $K$ of $\cU_q(\su )$, that is, 
\begin{equation}\label{podinv}
\cO(\Sq )=\{ x\in \cO(\SU)\,;\,x=x\anf K\}.
\end{equation}
From (\ref{tefkact}) and (\ref{podinv}), it follows that $\cO(\Sq )$   
is the linear span of elements $t^l_{0k}$, 
where $k=-l,{\dots},l$ and $l\in\dN_0$.

The importance of the Podle\'s sphere $\Sq $ stems from the fact that it 
is a quantum homogeneous space for the compact quantum group $\SU$. 
In fact, $\cO(\Sq )$  is a right $\ast$-coideal of the Hopf $\ast$-al\-ge\-bra 
$\cO(\SU)$  and hence a left $\cU_q(\su )$-module $\ast$-al\-ge\-bra 
with respect to the left action (\ref{act}). 
The generators
$$
x_{-1} := (1+q^{-2})^{1/2} B,\quad x_0 := 1-(1+q^2)A,\quad  
x_1 := -(1+q^2)^{1/2} B^\ast
$$
transform by the spin 1 matrix corepresentation of $\cO(\SU)$, see
\cite[p.\ 124]{KS} or \cite{P1}. We have
\begin{align}\label{delta1}
\Delta (B) &= B\otimes a^2 + (1-(1+q^2)A) \otimes ac -q B^\ast \otimes c^2,\\
\label{delta2}
\Delta (A) &= -q^{-2} B\otimes ab -q^{-1} (1-(1+q^2)A) \otimes bc 
+ A \otimes 1 -B^\ast \otimes dc,\\
\label{delta3}
\Delta (B^\ast) &= -q^{-1} B \otimes b^2 - (1-(1+q^2)A) \otimes db 
+ B^\ast \otimes d^2. 
\end{align}
The left crossed product algebra $\cO(\Sq )\rti \cU_q(\su )$ 
is the $\ast$-al\-ge\-bra generated by the 
$\ast$-sub\-al\-ge\-bras $\cO(\Sq )$  
and $\cU_q(\su )$ with respect to the cross relations 
$$
fx =(f_{(1)}  \ang x) f_{(2)}\equiv 
\langle f_{(1)}, x_{(2)}\rangle x_{(1)} f_{(2)},  
\quad f\in \cU_q(\su ),\ \, x\in\cO(\Sq ).
$$
From (\ref{pair}) and (\ref{delta1})--(\ref{delta2}), we obtain the 
following cross relations for the generators
\begin{align}                                         \label{cross1}
&KA = AK,\  EA = AE + q^{-1/2} B^\ast K,\  FA = AF - q^{-3/2} BK,\\
                                                     \label{cross2}
&KB = q^{-1} BK,\  EB = q BE + q^{1/2} (1-(1+q^2)A)K,\  FB = q BF,\\
                                                      \label{cross3}
&KB^\ast = qB^\ast K,\  EB^\ast = q^{-1} B^\ast E,\  
FB^\ast = q^{-1} B^\ast F - q^{-1/2} (1-(1+q^2)A)K.
\end{align}
For the Haar state $h$ on $\cO(\SU)$,  it is known that  
$h(xy) = h((K^{-2} \ang y \anf K^{-2}) x)$ for $x,y\in \cO(\SU)$. 
By a slight abuse of notation, we denote the restriction of $h$ to 
$\cO(\Sq )$  also by $h$. Then $h$ is the unique $\cU_q(\su )$-in\-var\-i\-ant 
state on the left $\cU_q(\su )$-module $\ast$-al\-ge\-bra $\cO(\Sq )$. 
Since $\cO(\Sq )$ is right $K$-in\-var\-i\-ant, we have
\begin{equation}\label{haut}
h(xy) = h(\sigma (y) x)\ \,\mbox{with}\ \sigma(y) =K^{-2} \ang y,\quad 
x,y\in \cO(\Sq ). 
\end{equation}
%
\section{Dirac operator and differential calculus on ${\bf \Sq }$}  
\label{S3}
%
First we introduce linear subspaces $V^+$, $V^-$ and $V$ of 
the $\ast$-al\-ge\-bra $\cO(\SU)$ which are left $\cO(\Sq )$-modules 
and left $\cU_q(\su )$-modules. Set
\begin{align*}
&V^+ = {\rm Lin} \{ v^{l+1/2}_{1/2,k}\,;\, 
k=-(l\!+\!1/2),{\dots}, l\!+\!1/2,\ l\in \dN_0\},\\
&V^- = {\rm Lin} \{ v^{l+1/2}_{-1/2,k}\,;\, 
k=-(l\!+\!1/2),{\dots}, l\!+\!1/2,\ l\in \dN_0\}
\end{align*}
and $V=V^+\oplus V^-$, where 
$v^{l+1/2}_{\pm 1/2,k}=[2l\!+\!1]_q q^{\pm 1/2}t^{l+1/2}_{\pm 1/2,k}$ 
and $t^{l+1/2}_{\pm 1/2,k}$ are the 
matrix elements from Section \ref{S2}.
For brevity, we shall frequently write $f v^\pm$ for $f\ang v^\pm$, 
where $v^\pm\in V^\pm$.
\begin{thl}  \mbox{ }                \label{L4/3}          \label{L3/1}
\begin{enumerate}
  \item[(i)] 
        $V^+ \,=\, \cO(\Sq ) c+ \cO(\Sq ) d
        \,=\, c\,\cO (\Sq ) + d\,\cO(\Sq ),\sn\\
        V^- \,=\,\cO (\Sq ) a + \cO(\Sq ) b
           \,=\, a\,\cO (\Sq ) + b\, \cO(\Sq )$.
  \item[(ii)] 
        $f\ang V^\pm\subseteq V^\pm$\ \, for $f\in \cU_q(\su )$.
 \item[(iii)]
  $R_E(\cO(\Sq ))\, \subseteq\, 
    \cO(\Sq )a^2 + \cO(\Sq )ab + \cO(\Sq )b^2\sn\\
    \phantom{R_E(\cO(\Sq ))} \quad
    \,=\,a^2\, \cO(\Sq )+ ab\, \cO(\Sq ) + b^2\, \cO(\Sq ), \sn\\
    R_F(\cO(\Sq ))\,\subseteq\, 
     \cO(\Sq )c^2 + \cO(\Sq )cd + \cO(\Sq )d^2\sn\\
    \phantom{R_F(\cO(\Sq ))} \quad 
    \,=\,c^2\, \cO(\Sq )+ cd\, \cO(\Sq ) + d^2\, \cO(\Sq ).$
 \item[(iv)] 
    $R_E (x) R_F (y) \in\cO(\Sq )$,\ \ $R_F (x) R_E (y) 
    \in\cO(\Sq )$\ \,for $x,y \in\cO (\Sq )$.
\end{enumerate}
\end{thl}
{\bf Proof.} (i) The above formulas for the matrix elements imply 
that $t^{l+1/2}_{1/2,k}\subseteq \cO(\Sq ) d$ and 
$t^{l+1/2}_{1/2,-k} \subseteq \cO(\Sq )c$ for $k>0$, 
so $V^+ \subseteq \cO(\Sq )c+\cO(\Sq )d$. Conversely, 
from the form of the Clebsch-Gordan coefficients 
(see e.g.\ \cite[Section 3.4]{KS}) and the definition of $t^{l+1/2}_{1/2,j}$, 
it follows that $t^l_{0,k} c \subseteq V^+$ 
and $t^l_{0,k} d\subseteq V^+$ for all $l\in\dN_0$. 
Thus, $V^+ = \cO(\Sq )c + \cO(\Sq )d$. The other relations are proved 
by a similar argumentation.

(ii) follows from (i) and Equations (\ref{leftact1}) and 
(\ref{leftact2}), or from the formulas (\ref{EFK}).

(iii) From the definition of the dual pairing (\ref{pair}) 
and the formulas (\ref{delta1})--(\ref{delta3}), we obtain for the 
right action of $E$ and $F$ on the generators $A$, $B$, $B^\ast$
\begin{align*}
R_E(B)&=-q^{-1/2} a^2, & R_E(B^\ast) &= q^{-3/2} b^2, & 
R_E(A)&= q^{-3/2} ba,\\
R_F(B)&=-q^{3/2} c^2, & R_F(B^\ast) &= q^{1/2} d^2, & 
R_F(A)&= q^{1/2} dc.
\end{align*}
Using Equation (\ref{actmod}) and the $K$-invariance (\ref{podinv}) 
of $\cO(\Sq )$,  
assertion (iii) follows easily by induction on powers of generators. 

(iv) follows at once from (iii) because pairwise products of the elements 
$a^2$, $ab$, $b^2$ with $c^2$, $cd$, $d^2$ are contained 
in $\cO(\Sq )$.                                      \hfill $\Box$
\sn

From Lemma \ref{L3/1}(i) and (ii), we conclude that $V^+$ and $V^-$ are 
invariant under the action of the {\it left} crossed product algebra 
$\cO(\Sq )\rti \cU_q(\su )$.
It can be shown that the corresponding modules are 
irreducible $\ast$-re\-pre\-sen\-tations of the cross product 
$\ast$-al\-ge\-bra $\cO(\Sq ) \rti \cU_q(\su )$.

Now we define the Dirac operator $D$ on $\Sq $. In order to do so, 
we essentially use the {\it right} action of $\cU_q(\su )$ on $\cO(\SU)$. 
Since $\alpha^{l+1/2}_{1/2} = [l+1]_q$, it follows from (\ref{tefkact}) that
\begin{equation}                                             \label{efeig}
R_E (v^{l+1/2}_{1/2,k}) = -[l\!+\!1]_q\,v^{l+1/2}_{-1/2,k}, \quad 
R_F(v^{l+1/2}_{-1/2,k}) = -[l\!+\!1]_q\, v^{l+1/2}_{1/2,k}
\end{equation}
for $k = -(l+\frac{1}{2}),{\dots},l+\frac{1}{2}$. 
Hence we have $R_E: V^+\rightarrow V^-$ and  
$R_F=(R_E)^\ast : V^-\rightarrow V^+$. 
Therefore, the operator $D$, given by the matrix 
\begin{equation*}
D := \left( \begin{matrix} 0 &R_F\\
R_E &0 \end{matrix} \right) 
\end{equation*}
on $V=V^+\oplus V^-$, maps $V$ into itself and is hermitian. 
Moreover, by (\ref{efeig}), 
$\psi^\epsilon_{l+1,k} := \frac{1}{\sqrt{2}} (v^{l+1/2}_{1/2,k}, 
\epsilon v^{l+1/2}_{-1/2,k} )^\TT $, 
$k=-(l+\frac{1}{2}),{\dots},l+\frac{1}{2},$ 
is an orthonormal sequence of eigenvectors of the operator $D$ 
with respect to the eigenvalues $-\epsilon [l+1]_q$, where $\epsilon = \pm1$ 
and $l\in\dN_0$. Let $\cK$ denote the closure of the subspace $V$ in the 
Hilbert space $\cL^2 (\SU)$. 
By the preceding, the closure of $D$,
denoted again by $D$, 
is a self-adjoint operator on $\cK$,
has a bounded inverse and $|D|^{-z}$ is 
of trace class for all $z\in \dC$, $\RT \, z>0$. 
The vectors $\varphi^+_{l+1,k} := (v^{l+1/2}_{1/2,k}, 0)^\TT $,  
$\varphi^-_{l+1,k} := (0,v^{l+1/2}_{-1/2,k})^\TT $, 
$k= -(l+\frac{1}{2}),{\dots},l+\frac{1}{2}$, $l\in\dN_0$,  
form an orthonormal basis of the Hilbert space $\cK$ consisting of 
eigenvectors of the self-adjoint operator $|D|=(D^\ast D)^{1/2}$ 
with respect to the eigenvalues $[l+1]_q$.

Let $x\in \cO(\Sq )$. 
Using (\ref{actmod}), (\ref{rightact2}) and (\ref{podinv}), we compute 
\begin{align*}
R_F(xv^-) &=-q xv^- \anf F
= -q (x\anf F)(v^-\anf K) -q(x\anf K^{-1}) (v^-\anf F)\\
&=q^{-1/2} R_F (x)  v^- +x R_F (v^-),
\end{align*}
so $[R_F, x] v^- = q^{-1/2} R_F(x) v^-$ for $v^- \in V^-$. 
Notice that $R_F(x) v^-\in V^+$ by Lemma \ref{L4/3} and 
Equations (\ref{t1})--(\ref{t4}). 
Similarly, we have for all $v^+\in V^+$
the identity $[R_E, x] v^+ = q^{1/2} R_E (x) v^+\in V^-$. 
Thus, for all $(v^+,v^-)^\TT \in V^+\oplus V^-$, we obtain  
$[D,x] (v^+,v^-)^\TT = (q^{-1/2} R_F (x) v^-, q^{1/2} R_E (x)v^+)^\TT$, 
that is,
\begin{equation}                                     \label{dcom}
[D,x] = \left( \begin{matrix} 0 &q^{-1/2} R_F(x)\\
q^{1/2} R_E (x) &0 \end{matrix}\right),
\end{equation}
where the elements $R_E(x)$ and $R_F(x)$ of $\cO(\SU)$ act by left 
multiplication on $V^-$ and $V^+$, respectively. In particular,  
$[D,x]$ is a bounded operator on the Hilbert space $\cK$ for 
all $x\in \cO(\Sq )$. Since $D$ is a self-adjoint operator with 
compact resolvent as noted above, $(D,\cO(\Sq ),\cK)$ is a spectral 
triple \cite{C,GVF}.

We turn now to the reality structure. 
Let $J$ denote the involution of the $\ast$-al\-ge\-bra $\cO(\SU)$. 
By Lemma \ref{L3/1}(i), $J$ maps $V^\pm$ onto $V^\mp$. More precisely, 
from the formulas for the elements  $t^{l+1/2}_{\pm 1/2,j}$ listed 
in Section \ref{S2} and  the fact that the polynomials 
$p_k(\cdot;\cdot,\cdot)$ have real 
coefficients it follows that
\begin{equation}                               \label{tinvol}
(t^{l+1/2}_{\pm 1/2,j})^\ast = (-q)^{j\mp 1/2} t^{l+1/2}_{\mp 1/2,-j}. 
\end{equation}
Define an anti-linear operator $\cJ_0:=\im K R_{K^{-1}}J$ on $V$, that is, 
$$
\cJ_0(v)=\im K\ang v^\ast \anf K,\quad v\in V.
$$ 
By (\ref{actstar}), 
$\cJ_0^2(v)=\im K\ang (\im K\ang v^\ast \anf K)^\ast \anf K
=K\ang (K^{-1} \ang v\anf K^{-1})\anf K=v$ for  $v\in V$, and so
$\cJ_0^2=I$. 
The same arguments imply 
\begin{equation}                                \label{Jcom}
\cJ_0 y^\ast\cJ_0^{-1}v= v (K\ang y \anf K),\quad 
y\in \cO(\Sq ),\ v\in V, 
\end{equation}
so that 
$x \cJ_0 y^\ast\cJ_0^{-1}v=\cJ_0 y^\ast\cJ_0^{-1}x v
= xv (K\ang y \anf K)$ 
for all $x,y \in\cO(\Sq )$. 
Hence $[x,\cJ_0 y^\ast\cJ_0^{-1}]=0$. Also, by (\ref{dcom}) and (\ref{Jcom}), 
$[[D,x],\cJ_0 y^\ast\cJ_0^{-1}]=0$ which means that $\cJ_0$ satisfies the 
order-one condition. 

From (\ref{tefkact}), (\ref{EFK}) and (\ref{tinvol}), we get
$$
\cJ_0(v^{l+1/2}_{\pm 1/2,j})
=\im [2l\!+\!1]_q^{1/2} q^{\pm 1/2}
(-q)^{j\mp 1/2} K \ang t^{l+1/2}_{\mp 1/2,-j}\anf K 
=\pm  \im ^{2j} v^{l+1/2}_{\mp 1/2,-j}
$$
for $j=-(l+1/2),{\dots},l+1/2$ and $l\in\dN_0$. 
Remind that $(-1)^{j\mp 1/2}$ is real and $2j$ is an odd integer. 
Since $\{v^{l+1/2}_{\pm 1/2,j}\}$ forms an orthonormal basis of the 
Hilbert space $\cK$, the operator $\cJ_0$ is anti-unitary on $\cK$. 
It is easy to check that $\cJ_0$ anti-commutes with the Dirac operator $D$.

Recall that the spectral triple $(D,\cO(\Sq ),\cK)$ 
is called {\it even} \cite[Definition 1]{C2}, if there 
is a grading operator $\gamma$ on $\cK$ such that
$\gamma^\ast=\gamma$, 
$\gamma^2=I$, $\gamma D= - D \gamma$ and 
$\gamma x=x\gamma $ for all $x\in\cO(\Sq )$. 
With $\gamma$ defined by  $\gamma  v^{\pm}=(\pm 1)v^\pm$, 
$v^\pm\in V^\pm$, the spectral triple 
$(D,\cO(\Sq ),\cK)$ is even. 
To obtain a real structure on $(D,\cO(\Sq ),\cK)$, 
we change $\cJ_0$ slightly and set $\cJ :=  \gamma \cJ_0$. 
Then $\cJ$ is an anti-unitary operator, 
${\cJ}^2 =-I$,  $\gamma \cJ=-\cJ \gamma$
and $\cJ D=D\cJ$.
Moreover,  $[x,\cJ y^\ast\cJ^{ -1}]=0$ 
and $[[D,x],\cJ y^\ast\cJ^{ -1}]=0$. 
This means that  $\cJ$ is a real structure on the 
even spectral triple $(D,\cO(\Sq ),\cK)$ 
according to \cite[Definition 3]{C2}.

Next we construct a first order differential calculus on $\Sq $. 
Let $\Gamma = \cL(V)$ be the algebra of linear operators on $V$. 
By Lemma \ref{L3/1}(i), $\Gamma$ is an $\cO(\Sq )$-bimodule 
with respect to the left and right multiplication of 
the algebra $\cO(\Sq )$. Defining $\dd:\cO(\Sq )\rightarrow \Gamma=\cL(V)$ 
by $\dd x := \im [D,x]$, we obtain a first order differential $\ast$-calculus 
$(\Gamma,\dd)$ on the algebra $\cO(\Sq )$. Since the right actions $R_E$ 
and $R_F$ commute with the left action of $f\in\cU_q(\su )$ 
on $\cO(\SU)$, the operator $D$ and so  
the differentiation $\dd$ commute with the left $\cU_q(\su )$-action. 
The latter means that 
the differential calculus $(\Gamma,\dd)$ is covariant with respect 
to the left $\cU_q(\su )$-action.

We show that $(\Gamma,\dd)$ has the quantum tangent space 
$\cT = {\rm Lin} \{E,F\}$. Indeed, 
let $x_j$, $y_j$, $j=1,\dots,n$, 
be such that $\sum_j \dd x_j\, y_j=0$. Then 
$$
\sum_j \dd x_j\, y_j\, (v^+, v^-)^\TT
=\sum_j \im (q^{-1/2} R_F (x_j)y_j v^-,q^{1/2} R_E (x_j) y_j v^+)^\TT=0
$$
for all $v^\pm \in V^\pm$ or, equivalently,
$\sum_j R_E(x_j) y_j = \sum_j R_F (x_j) y_j =0$. Hence 
$
\sum_j x_{j(1)}\otimes x_{j(2)} y_j\in\cT^\bot\otimes\cO(\SU),
$
where $\cT^\bot = \{x\in \cO(\Sq )\,;\,
\langle E, x\rangle=\langle F, x\rangle=0\}$. 
This shows that $\cT = {\rm Lin} \{E,F\}$ is the quantum tangent space 
of $(\Gamma,\dd)$. 
Since $\dim \cT=2$, $(\Gamma,\dd)$ is 2-dimensional.

The universal differential calculus $(\DCS ^\wedge_u ,\dif )$ 
associated with $(\Gamma,\dd)$  is defined by
$\DCS ^\wedge_u =\bigoplus _{k=0}^\infty \DCS ^{\ot k}/\mathcal{I}$, where
$\DCS ^{\ot k}=\DCS \otp \DCS \otp \cdots \otp \DCS $
($k$-fold tensor product), $\DCS ^{\ot 0}=\podl $ and 
$\mathcal{I}$ denotes the 2-sided ideal of the tensor algebra 
$\bigoplus _{k=0}^\infty \DCS ^{\ot k}$
generated by the subset
$\{\sum _i\dif x_i\ot \dif y_i\,;\,\sum _i\dif x_i\,y_i=0\}$. 
The product of the algebra $\DCS ^\wedge $ is
denoted by $\wedge $.
Let $(\Gamma^\wedge,\dd)$ be  the 
higher order calculus on $\cO(\Sq )$
obtained from $(\DCS ^\wedge_u ,\dif )$ 
by setting $\Gamma^{\wedge k}=0$ for $k\geq 3$. 
(It can be shown that $\DCS ^{\wedge k}_u=0$ for $k\geq 3$ but 
we do not need this here.) A covariant differential calculus 
$(\Gamma^\wedge_P,\dd_P)$ on $\cO(\Sq )$
was constructed by P. Podle\'s in \cite{P2}.
The following lemma is proved in the appendix. 
\begin{thl}                                    \label{app}
\begin{enumerate}
\item[(i)]  $(\Gamma,\dd)$ is the unique 2-dimensional covariant
first order differential calculus on $\cO(\Sq )$
such that $\{\dd A,\dd B,\dd B^\ast\}$ generate 
the right (left) $\cO(\Sq )$-module $\Gamma$.
\item[(ii)] 
The differential calculi $(\Gamma^\wedge,\dd)$ and $(\Gamma^\wedge_P,\dd_P)$ 
are isomorphic. 
\end{enumerate}
\end{thl}

The next theorem collects some properties of the Dirac operator $D$.
\begin{tht}
\begin{enumerate}
\item[(i)]  $D$ is a self-adjoint operator with discrete spectrum 
consisting of eigenvalues $\epsilon [n]_q$ with multiplicities $2n$, 
where $\epsilon =\pm 1$ and $n\in\dN$. 
In particular, $|D|^{-z}$ is of trace class for $z\in\dC$, $\RT\, z>0$.
\item[(ii)] 
 $D$ commutes with the left action of $\cU_q(\su )$ on $V$.
\item[(iii)] 
Equipped with the grading operator $\gamma$, 
the spectral triple $(D,\cO(\Sq ),\cK)$ is even. 
The operator $\cJ$ defines a real structure on $(D,\cO(\Sq ),\cK)$
according to \cite[Definition 3]{C2}. 
\item[(iv)] 
There is a unique 2-dimensional covariant
first order differential calculus $(\Gamma,\dd)$
on $\cO(\Sq )$ such that $\{\dd A,\dd B,\dd B^\ast\}$ generate 
the left $\cO(\Sq )$-module $\Gamma$. This calculus is given by 
given by $\dd x := \im [D,x]$, $x\in\cO(\Sq )$.
\end{enumerate}
\end{tht}

We close this section by giving a purely algebraic construction of the 
diffe\-ren\-tial calculus $(\Gamma,\dd)$. 
Since the subset $\cS=\{ b^k c^n\,;\,k,n\in\dN_0\}$ is a left and right 
Ore set and the algebra  $\cO (\SU)$ has no zero divisors, the 
localization algebra $\hat{\cO} (\SU)$ of $\cO(\SU)$ at $\cS$ exists. 
Then $\cO(\SU)$ is a $\ast$-sub\-al\-ge\-bra of $\hat{\cO}(\SU)$ and 
$b$ and $c$ are invertible in the larger algebra $\hat{\cO}(\SU)$. 
The crucial observation is Equation (\ref{efcom}) in the following lemma.
\begin{thl} \label{L3/2}
\begin{enumerate}
\item[(i)] 
The mappings $ R_F, R_E :\cO(\Sq )\rightarrow \cO(\SU)$ are 
derivations, that is, $R_f (xy)= x R_f (y) + R_f (x) y$ for $f=F,E$ and 
$x,y\in\cO(\Sq )$.
\item[(ii)]
For all $x\in \cO(\Sq )$,
\begin{equation}                                           \label{efcom}
R_F(x) =[q^{1/2} \lambda^{-1} db^{-1},x],\quad  
R_E (x) = -[q^{-1/2} \lambda^{-1} ac^{-1},x],
\end{equation}
where the commutators are taken in the algebra $\hat{\cO}(\SU)$.
\end{enumerate}
\end{thl}
{\bf Proof.} (i) follows from (\ref{actmod}) and 
the comultiplications of $F$ and $E$ 
combined with the $K$-invariance (\ref{podinv}) of $\cO(\Sq )$.

(ii) By (i), $R_F$ and $R_E$ are derivations. The commutators on the 
right hand side of (\ref{efcom}) are also derivations. Hence it suffices to 
check that Equation (\ref{efcom}) holds for the generators $x=A,B,B^\ast$. 
We omit the details of this easy computation.       \hfill $\Box$
\sn

From relations (\ref{dcom}) and (\ref{efcom}), it follows that, 
for all $x\in\cO(\Sq )$,
\begin{equation*}
\dd x = \im [D,x] = \im \lambda^{-1} \left( \begin{matrix} 0 &[db^{-1},x]\\
-[ac^{-1},x] &0 \end{matrix} \right).
\end{equation*}

The reason for the somehow surprising identity (\ref{efcom}) is the 
following fact \cite{SW3}: The elements  
$q^{1/2}\lambda K^{-1} E+K^{-2} c^{-1} a$ and 
$q^{1/2} \lambda FK^{-1} - qdb^{-1} K^{-2}$ of the right crossed 
product algebra $\cU_q(\su )\lti\hat\cO(\SU)$ commute with all elements 
of $\cO(\SU)$. By the $K$-invariance (\ref{podinv}), this 
implies (\ref{efcom}).

%
\section{The twisted cyclic cocycle on ${\bf \Sq }$} \label{S4}
We begin with a couple of preliminary lemmas.
\begin{thl}  \label{L4/4}
For all $x,y\in\cO(\Sq )$, we have
\begin{equation}  \label{heffe}
h(R_F(x) R_E (y)) = q^2 h(R_E(x) R_F(y)).
\end{equation}
\end{thl}
{\bf Proof.}
Let $x,y\in \cO(\Sq )$. 
By (\ref{actstar}) and (\ref{R}), 
$R_F(x)^\ast=x^\ast\anf E =-qR_E(x^\ast)$. 
The $K$-invariance (\ref{podinv}) of $\cO(\Sq )$ 
implies $(R_{F}R_{E}- R_{E}R_{F}) (y)=0$. Hence\\
\parbox{0.4cm}{\begin{align*}
 &&\\
 &&\\
 &&
\end{align*}}\hfill
\parbox{11cm}{
\begin{align*} 
h(R_F(x) R_E (y))&=\big(R_E (y), R_F(x)^\ast  \big) 
=-q\big(R_E (y), R_E(x^\ast ) \big)\\ 
&=-q\big(R_F R_E (y), x^\ast \big)
=-q\big(R_F (y), R_F (x^\ast) \big) \\
&=q^2\big(R_F (y), R_E (x)^\ast \big) 
=q^2h(R_E(x) R_F (y)).
\end{align*} }\hfill
\parbox{0.4cm}{\begin{align*}
 &\\
 &\\
 &\Box
\end{align*}}

For $x_0, x_1, x_2 \in \cO (\Sq )$, 
we define our cocycle $\tau$ by 
\begin{equation}                                 \label{cocdef}
\tau (x_0, x_1,x_2) 
= h(x_0 (R_F (x_1) R_E (x_2) - q^2 R_E (x_1) R_F (x_2))).
\end{equation}
Note that the products $R_F (x_1) R_E (x_2)$ and $R_E (x_1) R_F (x_2)$ are in 
$\cO(\Sq )$ by Lemma \ref{L4/3}(iv).
\begin{thl}                                             \label{L4/5}
Let $\sigma$ denote the algebra automorphism of $\cO(\Sq )$ given by 
$\sigma (x) = K^{-2} \ang x$, $x\in\cO(\Sq )$. 
Then $\tau$ is a non-trivial $\sigma$-twisted cyclic 
2-cocycle on $\cO(\Sq )$.
\end{thl}
{\bf Proof.} First, consider 
$\tau_1 (x_0,x_1,x_2) := h(x_0 R_F (x_1) R_E (x_2))$ and 
$\tau_2 (x_0, x_1,x_2) := h (x_0 R_E (x_1) R_F (x_2))$. 
Recall that $R_F$ and $R_E$ act as 
derivations on $\cO(\Sq )$ by Lemma \ref{L3/2}, 
and $h$ satisfies condition (\ref{haut}). 
Applying the boundary operator $b_\sigma$ to $\tau_1$ and $\tau_2$ 
and using the Leibniz rule for the derivations $R_E$ and $R_F$, 
one sees that the sum in (\ref{b1}) telescopes to zero. 
Hence $b_\sigma \tau=0$. 
Next, 
\begin{align*}
\tau (x_0,x_1,x_2)&=h (R_F (x_0x_1) R_E (x_2) - q^2 R_E (x_0 x_1) R_F (x_2))\\
& -h (R_F (x_0) R_E (x_1 x_2) - q^2 R_E (x_0) R_F (x_1 x_2))\\
& +h(R_F (x_0) R_E (x_1)x_2 - q^2 R_E(x_0) R_F (x_1)x_2)
=\tau (\sigma(x_2),x_0,x_1),
\end{align*}
where the first equality is an application of the Leibniz rule 
and the second equation follows from Lemma \ref{L4/4} and 
condition (\ref{haut}). 
Thus $\tau = \lambda_\sigma \tau$. 

We prove that $\tau$ is non-trivial. Consider the element 
\begin{align*}
\eta &:= B^\ast \otimes A \otimes B + q^2 B\otimes  B^\ast \otimes A 
+ q^2 A \otimes B \otimes B^\ast -q^{-2} B^\ast \otimes B \otimes A \\
&\quad - q^{-2} A\otimes  B^\ast \otimes B 
- B \otimes A \otimes B^\ast
+(q^6 - q^{-2}) A \otimes A \otimes A
\end{align*}
of the tensor product $\cO(\Sq )^{\otimes 3}$. 
By (\ref{leftact2}) and the definition of $\sigma$, we have
$\sigma(A)=A$, $\sigma(B)=q^2B$ and $\sigma(B^\ast)=q^{-2}B^\ast$. Since 
$\tau = \lambda_\sigma \tau$, it follows that 
$$
\tau (\eta) = 3\tau (B^\ast, A,B) - 3q^{-2}\tau (B^\ast,B,A)
+(q^6-q^{-2})\tau(A,A,A).
$$
From the definition of $\tau$, we obtain 
\begin{align*}
\tau (B^\ast, A, B) &= (q^2 - q^{-4}) h(A^3 - A^2) + q^{-2} h(A^2-A),\\
\tau (B^\ast, B, A) &= (q^4 - q^{-2}) h(A^3 - A^2) - q^{2} h(A^2-A),\\
\tau (A, A, A) &= (q^{-2} - q^4) h(A^3)- (q^{-2}-q^2) h(A^2).
\end{align*}
Inserting the values $h(A^j)\!=\!(1\!-\!q^2)/(1\!-\!q^{2j+2})$ 
and summing up gives 
$\tau(\eta)\!=\!-1$.

On the other hand, one computes 
$b_\sigma(\eta)= 2(q^4{-} q^{-2}) A\otimes A$ by 
using algebra relations (\ref{algrel}). 
Note that Equation (\ref{lambda}) and $\sigma(A)=A$ imply 
$\tau^\prime (A,A) =0$ for any 
$\sigma$-twisted cyclic 1-cocycle  $\tau^\prime$. 
If there 
were a $\sigma$-twisted cyclic 1-cocycle  $\tau^\prime$ such that 
$\tau =b_\sigma (\tau^\prime)$, we would get 
$$
\tau(\eta) = (b_\sigma \tau^\prime)(\eta)=\tau^\prime (b_\sigma \eta)
=2(q^4 {-} q^{-2}) \tau^\prime (A,A) =0,
$$
a contradiction. Thus $\tau$ is non-trivial.                \hfill $\Box$
\sn

The left action of $f\in\cU_q(\su )$ on cycles 
$\eta = \sum_k x^k_0\otimes {\cdots}\otimes x^k_n 
\in \cO(\Sq )^{\otimes n+1}$ is defined by 
$$
  f\ang \eta 
=\sum_k f_{(1)}\ang x^k_0\otimes {\cdots}\otimes f_{(n+1)}\ang x^k_n.
$$
Pairing a 2-cycle $\eta$  
with $\tau$ gives $\tau(f\ang \eta)=\varepsilon(f)\tau(\eta)$ 
since the  right actions $R_E$ 
and $R_F$ commute with the left action of $\cU_q(\su )$ on 
$\cO(\SU)$ and $h$ is $\cU_q(\su )$-in\-var\-i\-ant. 
Hence $\tau$ is $\cU_q(\su )$-in\-var\-i\-ant. 

We next describe $\tau$ analytically. 
Let $\cK^\pm$ be the closure of $V^\pm$ in the Hilbert space 
$\cL^2(\SU)$ and, for $z\in \dC$, $\RT\, z>2$, 
let $\zeta(z)$ denote the holomorphic function given by
$$
\zeta(z) = \sum^\infty_{n=1} [n]_q^{-z} [2n]_q.
$$
The following lemma is a slight modification of Theorem 5.7 in \cite{SW2}.
\begin{thl}                                                \label{L4/6}
Let  $z\in\dC$ and $\RT\, z>2$. For any $x\in \cO(\Sq )$, 
the closure of the operator $K^2 |D|^{-z} x$ restricted to 
the Hilbert spaces $\cK^\pm$ is of trace class and 
\begin{equation}\label{htrace}
h(x) = \zeta (z)^{-1} \tr_{\cK^\pm} K^2|D|^{-z} x.
\end{equation}
\end{thl}
{\bf Proof.} The operators $|D|^{-z}$ and $K$ act on the orthonormal basis 
$\{ v^{l+1/2}_{\pm 1/2,k}\}$ of the Hilbert spaces $\cK^\pm$ by 
$|D|^{-z} v^{l+1/2}_{\pm 1/2,k} = [l+1]_q^{-z} v^{l+1/2}_{\pm 1/2,k}$ and 
$K v^{l+1/2}_{\pm 1/2,k}=q^k v^{l+1/2}_{\pm 1/2,k}$, respectively. 
Since $x\in \cO(\Sq )$ acts as a bounded operator on 
$\cK^\pm$, $K^2|D|^{-z} x$ is of trace class. 
Thus, $h_z (x) := \tr_{\cK^\pm}  K^2 |D|^{-z} x$ 
is well defined for $x\in\cO(\Sq )$. 

We show that $h_{2z}$ is a $\cU_q(\su )$-in\-var\-i\-ant linear functional on 
$\cO(\Sq )$, that is, $h_{2z} (f\ang x) = \varepsilon(f)h_{2z} (x)$ for 
$f\in\cU_q(\su )$ and $x\in\cO(\Sq )$. In order to do so, we essentially 
use the fact that $V^\pm$ are left modules of the left crossed product 
algebra $\cO(\Sq ) \rti\cU_q(\su )$. It suffices to verify the invariance 
for the generators $f=E,F, K,K^{-1}$. We carry out the proof for $f{=}E$ 
and show that $h_{2z} (E\ang x){=}0$.

From the relations (\ref{cross1})--(\ref{cross3}) of the cross product 
algebra $\cO(\Sq )\rti \cU_q(\su )$, it follows that for any $x\in\cO(\Sq )$ 
there are elements $x,x^\prime, x^{\prime\prime}\in \cO(\Sq )$ 
such that $Kx=yK$ and $xE=Ex^\prime + Kx^{\prime\prime}$. Thus, 
since $E v^{l+1/2}_{\pm 1/2,k} = \alpha^{l+1/2}_k v^{l+1/2}_{\pm 1/2,k+1}$ 
and $|\alpha^{l+1/2}_k|\le \mbox{const}\, q^l$, 
the operators $K^2 |D|^{-2z} ExK^{-1}$, $|D|^{-z} KxE$ and $KE|D|^{-z}$ are 
of trace class. Since $|D|^{-z}$ commutes with $E$ and $K$, we compute
\begin{align}
&\tr_{\cK^\pm}  (K^2 |D|^{-2z} E x K^{-1}) 
\!=\! \tr_{\cK^\pm} (|D|^{-2z} K^2 EK^{-1} y)     
\!=\! q \tr_{\cK^\pm}  (KE|D|^{-z}) |D|^{-z} y         \nonumber \\
& \!=\! q \tr_{\cK^\pm}  |D|^{-z} y (KE |D|^{-z}) 
\!=\! q \tr_{\cK^\pm}  (|D|^{-z} y KE) |D|^{-z} 
\!=\! q \tr_{\cK^\pm}  (|D|^{-2z} y KE)       \nonumber \\
& \!=\! q \tr_{\cK^\pm}  (|D|^{-2z} K x E) 
\!=\! q \tr_{\cK^\pm}(K^2 |D|^{-2z} K^{-1} xE),                \label{kerel}
\end{align}
where all operators in parentheses are of trace class.
(Strictly speaking, one has to take the closures of these operators.)
In the left crossed product algebra we have $f\ang x= f_{(1)} x S(f_{(2)})$ 
and so $E\ang x = E x K^{-1} -q K^{-1} x E$.
Hence $h_{2z} (E\ang x)=0$ by (\ref{kerel}).

The functional $h_{2z}$ on $\cO(\Sq )$ is $\cU_q(\su )$-in\-var\-i\-ant. Since 
\begin{equation*}
h_z(1) = \tr_{\cK^\pm} K^2 |D|^{-z} 
\!=\! \sum^\infty_{l=0}\,\sum^{l+1/2}_{k=-(l+1/2)} [l+1]_q^{-z} q^{2k} 
\!=\! \sum^\infty_{l=0} [l+1]_q^{-z} [2l+2]_q=\zeta (z),
\end{equation*}
it follows that $\zeta (2z)^{-1} h_{2z}$ is the invariant state $h$ 
on $\cO(\Sq )$. Finally, $\zeta (z) h(x)$ and $h_z (x)$ are 
holomorphic functions for $z\in\dC$, $\RT\, z>2$. Since they are equal 
for $\RT\, z>4$ as just shown, they coincide also for $\RT\, z>2$.
                                                            \hfill $\Box$
\sn

Using Lemma \ref{L4/6}, we now express the cocycle $\tau$ defined by 
(\ref{cocdef}) in terms of our Dirac operator $D$. Let
\begin{equation*}
\gamma_q = \left( \begin{matrix} 1 & 0\\
0 &-q^2 \end{matrix}\right)
\end{equation*}
be the ``grading'' operator on $\cK = \cK^+\oplus \cK^-$. 
Indeed, using formulas (\ref{dcom}), (\ref{cocdef}) and (\ref{htrace}),  
we obtain
\begin{align}
&\tr_{\cK}\,\gamma_q K^2 |D|^{-z} x_0 [D, x_1][D, x_2]\nonumber\\
&\quad =  \tr_{\cK^+}  
K^2 |D|^{-z} x_0 R_F (x_1) R_E (x_2) 
- q^2 \tr_{\cK^-} K^2 |D|^{-z} x_0 R_E (x_1) R_F (x_2)\nonumber\\
&\quad 
=\zeta (z)\, h(x_0 R_F (x_1) R_E (x_2) {-} q^2 x_0  R_E (x_1) R_F (x_2))
= \zeta (z)\, \tau(x_0,x_1,x_2)                        \label{res}
\end{align}
for $z\in\dC$, $\RT\, z>2$, and $x_0,x_1,x_2\in\cO(\Sq )$.
Using the binomial series, one can show that 
$$
\zeta (z)= (q^{-1}-q)^{z-1}\sum_{k=0}^\infty 
\binom{1-z}{k}\biggl((1-q^{2(z-2+k)})^{-1}+(1-q^{2(z-1+k)})^{-1}\biggr).
$$ 
The right-hand side is a meromorphic function which is denoted again by 
$\zeta(z)$. It has a simple pole at $z=2$ 
with residue $\lambda (\log q)^{-1}$. Therefore, by (\ref{res}), 
$$
\underset{z=2}{{\rm res}}~ \tr_{\cK}\, 
 \gamma_q K^2 |D|^{-z} x_0 [D,x_1][D,x_2]
=\lambda (\log q)^{-1} \tau (x_0,x_1,x_2).
$$

On the other hand, the cocycle $\tau$ can also be obtained from the 
differential calculus $(\Gamma^\wedge,\dd)\cong(\Gamma^\wedge_P,\dd_P)$. 
As shown in \cite{P2}, there is a 
left-in\-var\-i\-ant 2-form $\omega\ne 0$ such that 
$\Gamma^{\wedge 2} = \cO(\Sq )\omega$ and $x\omega=\omega x$ 
for $x\in\cO(\Sq )$. 
The following remarkable result is due to I. Heckenberger. 
His proof will be given in the appendix. 
\begin{thl}                                        \label{L4/7}
The volume form $\omega$ can be chosen such that 
$\tau$ is equal to the $\sigma$-twisted cyclic cocycle $\tau_{\omega,h}$.
\end{thl}

Retaining the preceding notation, we now summarize the main results 
obtained in this section.
\begin{tht} \label{T4/8}
There is a non-trivial $\sigma$-twisted cyclic 2-cocycle $\tau$ on the 
algebra $\cO(\Sq )$ such that 
\begin{align*}
\tau (x_0, x_1, x_2) &= h(x_0 (R_F (x_1) R_E (x_2) - q^2 R_E (x_1) R_F(x_2)))\\
&=\zeta(z)^{-1} \tr_\cK\, \gamma_q K^2 |D|^{-z} x_0 [D, x_1][D,x_2]\\
&=\lambda^{-1} (\log q)~ \underset{z=2}{{\rm res}}~ 
\tr_\cK\, \gamma_q K^2 |D|^{-z} \,x_0[D,x_1][D,x_2] 
\end{align*}
for $z\in\dC$, $\RT\, z>2$, and $x_0, x_1, x_2\in\cO(\Sq )$. 
The cocycle $\tau$ is $\cU_q(\su )$-in\-var\-i\-ant and 
coincides with the $\sigma$-twisted cyclic 
cocycle $\tau_{\omega,h}$ associated with the volume form $\omega$ 
of the differential calculus $(\Gamma, \dd)$.
\end{tht}

In "ordinary" non-commutative geometry the grading operator anti-commutes 
with the Dirac operator $D$ and with the anti-unitary operator 
implementing the real structure. 
This is not true for our grading operator $\gamma_q$,
but $\gamma=q^{-1} \gamma_qR_{K^{-2}}$ does  
anti-commute with $D$ and  $\cJ$.

%
%
\section{Appendix} \label{S5}

In this appendix we present the proofs of 
Lemmas \ref{app} and \ref{L4/7} by I. Heckenberger.
Let us first state some general facts on covariant differential calculi 
(see e.g.\ \cite{H}). 
To an (arbitrary) covariant first order differential 
calculus $(\DCS ,\dif )$ over
$\podl $, one associates the left ideal
$\lid :=\{b\in \podl \,;\,\vep (b)=\pinv (b)=0\}$, where the map 
$\pinv :\podl \to \DCS \otp \OSUq 2$ is given by
$\pinv (b)=\dif b_{(1)}\ot S(b_{(2)})$. 
The left ideal $\lid$ determines $(\DCS ,\dif )$ uniquely. 
The (right) quantum tangent space
associated with $\DCS $ is the linear subspace 
$\cT:=\{f\in \podl '\,;\,\pair{f}{1}=\pair{f}{\lid }=0\}$ of the  
dual vector space $\podl '$ of $\podl $. 
Set $\podl ^+=\{x\in \podl \,;\,\vep (x)=0\}$. 
The cardinal number  $\dim \DCS :=\dim _\comp \DCS /\DCS \podl ^+
=\dim _\comp \podl^+ /\lid$ is called the (right) dimension of $\DCS $. 
If $\dim \DCS <\infty $, then the quantum tangent space $\cT$ also 
determines $(\DCS ,\dif )$ uniquely 
and $\dim \DCS=\dim_\comp \cT$. 

From now on, $(\DCS ,\dif )$ stands for the first order differential 
calculus over $\cO(\Sq )$ with quantum tangent space 
$\cT=\mathrm{Lin}\{E,F\}$ constructed in Section \ref{S3}. 
\medskip 

\noindent
{\bf Proof of Lemma \ref{app}.}  (i) 
Let $(\tilde\DCS ,\dif )$ be a finite dimensional covariant 
first order differential 
calculus over $\cO(\Sq )$ with quantum tangent space $\tilde\cT$. 
By the right-handed version of \cite[Corollary 5]{HK3}, we have 
$\tilde\cT\comp[K,K^{-1}]\subset \tilde\cT$. 
Since $\dim \tilde\DCS<\infty$, we deduce 
from \cite[Theorem 6.5.1]{HK2} that 
$\tilde\cT\subset\cU_q(\mathrm{sl}_2)|_{\podl}$. 
Using once more \cite[Corollary 5]{HK3}, one checks that only the 
linear spaces $\Lin\{E,F\}$, $\Lin\{E,E^2\}$ and $\Lin\{F,F^2\}$ 
are quantum tangent spaces of 
2-dimensional covariant first order differential
calculi over $\podl $. 
However, only one of them, namely $\cT=\mathrm{Lin}\{E,F\}$,  
is separated by the right $\OSUq 2$-co\-mo\-dule 
$\mathrm{Lin}\{B,B^*,$ $1-(1+q^2)A\}$. 
Hence, by the right-handed version of \cite[Lemma 7]{HK3},     
$(\DCS ,\dif )$ is the unique 2-dimensional covariant 
first order differential calculus generated by 
$\cg:=\mathrm{Lin}\{\dif B,\dif B^*,\dif A\}$ as a right $\podl$-module.

(ii)  
We first show that the first order
differential calculi $(\DCS _P,\dif _P)$ and $(\DCS ,\dif )$ 
are isomorphic. 
By the definition in \cite{P2}, the right $\podl $-module $\DCS _P$ 
is also generated by $\cg$. 
Since $\DCS _P$
is the quotient of $\mathrm{Lin}\,\cg\ot \podl $ by a submodule generated 
by a non-zero element \cite{P2}, we deduce 
$\dim \DCS _P=2$. 
Let $\cT_P$ denote the right quantum tangent space of $(\DCS _P,\dif _P)$.
Again by \cite[Lemma 7]{HK3}, $\mathrm{Lin}\{B,B^*,1-(1+q^2)A\}$
separates $\cT_P$. 
Applying the uniqueness result 
from the proof of Lemma \ref{app}(i) 
shows that $(\DCS _P,\dif _P)$ and $(\DCS ,\dif )$ 
are isomorphic.

By definition, $\DCS^{\wedge k}=\DCS^{\wedge k}_P=0$ for $k\geq3$. 
To complete the proof of 
Lemma \ref{app}(ii), it remains to verify that 
$\DCS ^{\wedge 2}_P$ and $\DCS ^{\wedge 2}$ are isomorphic as 
right $\podl $-mo\-du\-les and 
right $\OSUq 2$-comodules. 
Since $\DCS ^{\wedge 2}$ is obtained 
from the universal differential calculus, 
$\DCS ^{\wedge 2}_P$ is a quotient of
$\DCS ^{\wedge 2}$. 
From the definition of $\DCS ^{\wedge 2}_P$ in \cite{P2}, 
it follows that 
$\dim _\comp \DCS^{\wedge 2}_P /\DCS^{\wedge 2}_P \podl ^+=1$.
Thus it suffices to show that
$\dim _\comp \DCS^{\wedge 2} /\DCS^{\wedge 2} \podl ^+=1$ which 
follows from the right-handed version of \cite[Proposition 3.11(iv)]{HK4}.
                                                \hfill$\Box$
\medskip

\noindent
{\bf Proof of Lemma \ref{L4/7}.} 
To apply the results of \cite{HK4}, the pairing 
in \cite[Subsection 2.3.5]{HK4} has to be replaced by 
\begin{align}                     
      & \hspace{45pt} \tpair{\cdot }{\cdot }:
(\cT\ot \cT_0)\times (\DCS \otp \DCS \otp \OSUq 2) \to \comp,  \nonumber\\ 
&\tpair{t\ot s}{\dif x \ot \dif y \ot z}:=\pair{t\ot s}{x_{(1)}\ot x_{(2)}y^+}
\vep (z)=\pair{ts_{(1)}\ot s^+_{(2)}}{x\ot y}\vep(z), \label{eq-tpairdef}
\end{align}
where $t\in \cT$, $s\in \cT_0:=\mathrm{Lin}\{K^{-1}E,K^{-1}F\}$, 
$x,y\in \podl $, 
$z\in \OSUq 2$ and $f^+=f-\vep (f)$.
Set $\cT_2:=\{\sum _it_i\ot s_i\in \cT\ot \cT_0\,;\,
\sum _it_is_i\in \cT\}$.
Similarly to \cite[Corollary 2.10]{HK4}, one shows that (\ref{eq-tpairdef}) 
induces a non-degenerate pairing 
$$
\tpair{\cdot }{\cdot }:
\cT_2\times \DCS ^{\wedge 2}/\DCS ^{\wedge 2}\podl ^+\to \comp.
$$
Observe that $\cT_2=\comp t_2$ with
$t_2:=q^2F\ot K^{-1}E-E\ot K^{-1}F$. 
By (\ref{eq-tpairdef}) and the 
right $K$-invariance (\ref{podinv}) of $\podl$, 
\begin{align} \notag
\tpair{t_2}{\dif x\wedge \dif y}
&=\pair{q^2FK^{-2}\ot K^{-1}E-EK^{-2}\ot K^{-1}F}{x\ot y}\\
&=\pair{F\ot E-q^2E\ot F}{x\ot y}\label{eq-tpair}
\end{align}
for all $x,y\in\podl$. 
For notational convenience, we rename the generators 
$a$, $b$, $c$ and $d$ of $ \OSUq 2$ by 
$u_{11}$, $u_{12}$, $u_{21}$ and $u_{22}$, respectively.   
Set 
\begin{align*}
(p_{ij})_{i,j=1,2}\hspace{-1pt}:=\hspace{-1pt}(S(u^i_2)u^2_j)_{i,j=1,2}
\hspace{-1pt}=\hspace{-1pt}
\begin{pmatrix}A & B^*\\ B & 1\!-\!q^2A\end{pmatrix},
\quad \omega \hspace{-1pt}:=\hspace{-1pt}
\sum _{i,j,k}q^{2-2i}\dif p_{ij}\wedge \dif p_{jk}p_{ki}.
\end{align*}
As $\Delta(p_{ij})=\sum _{k,l}p_{kl}\ot S(u^i_k)u^l_j$ and 
the calculus $(\DCS ,\dif )$ is covariant, 
the 2-form 
$\omega$ is right-coinvariant.
Using Equation (\ref{eq-tpair}) and 
$\vep (p_{ij})=\delta _{i2}\delta _{j2}$, one computes 
$$
\tpair{t_2}{\omega }\!=\!
\sum _j\tpair{t_2}{q^{-2}\dif p_{2j}\wedge \dif p_{j2}}
=\sum _j\pair{F\ot E-q^2E\ot F}{q^{-2}p_{2j}\ot p_{j2}}\!=\!1.
$$
In particular,  $\omega \not=0$. 
By the definition in \cite{P2}, 
the right $\podl $-module $\DCS _P^{\wedge 2}\cong \DCS^{\wedge 2}$
is free and  generated by a non-zero right-coin\-var\-i\-ant 
central element.
Consequently,  $\DCS ^{\wedge 2}=\omega \podl $. 
Moreover, $\tpair{t_2}{\omega u_{(1)}}u_{(2)}
=\tpair{t_2}{\omega}\vep(u_{(1)})u_{(2)}=u$ 
for all $u\in \podl $. 
Let $\rho=\omega u\in \Gamma^{\wedge 2}$. Then, by the last relation and the 
coinvariance of $\omega$, 
$\rho =\omega \tpair{t_2}{\omega u_{(1)}}u_{(2)}
=\omega \tpair{t_2}{\rho _{(1)}}\rho _{(2)}$. 
This identity  and Equation (\ref{eq-tpair}) imply
\begin{align*}
x\dif y \wedge \dif z
&=x\omega \tpair{t_2}{\dif y_{(1)}\wedge \dif z_{(1)}}y_{(2)}z_{(2)}\\
&=x\omega \pair{F\ot E-q^2E\ot F}{y_{(1)}\ot z_{(1)}}y_{(2)}z_{(2)}\\
&=x\big(R_F(y)R_E(z)-q^2R_E(y)R_F(z)\big)\omega ,
\end{align*}
which proves Lemma \ref{L4/7}.
                                                      \hfill$\Box$

\end{document}